\documentclass[12pt,a4paper,oneside]{article}
\usepackage[english]{babel}
\usepackage[utf8]{inputenc}
\usepackage{amsmath}
\usepackage{amsfonts}
\usepackage{amssymb}
\usepackage{makeidx}
\usepackage{graphicx}
\usepackage{fourier}
\usepackage[left=2cm,right=2cm,top=2cm,bottom=2cm]{geometry}
\usepackage{authblk}

\numberwithin{equation}{section}

\newcommand{\ba}{\begin{array}\numberwithin{equation}{section}}
\newcommand{\ea}{\end{array}}
\newcommand{\bt}{\begin{tabular}}
\newcommand{\et}{\end{tabular}}
\newcommand{\btb}{\begin{table}}
\newcommand{\etb}{\end{table}}
\newcommand{\bc}{\begin{center}}
\newcommand{\ec}{\end{center}}
\newcommand{\bea}{\begin{eqnarray}}
\newcommand{\eea}{\end{eqnarray}}
\newcommand{\Bea}{\begin{eqnarray*}}
\newcommand{\Eea}{\end{eqnarray*}}

\newcommand{\beq}{\begin{equation}}
\newcommand{\eeq}{\end{equation}}

\title{Why the Kemeny Time is a Constant}
\author[1]{Karl Gustafson}
\author[2]{Jeffrey J. Hunter}
\date{3 November 2015}

\affil[1]{Department of Mathematics, University of Colorado, Boulder, Colorado, 80309-0395, USA.
 Email: karl.gustafson@colorado.edu}

\affil[2]{School of Engineering, Computer and Mathematical Sciences, Auckland University of Technology, New Zealand.
 Email: jeffrey.hunter@aut.ac.nz}

\begin{document}
\maketitle

\begin{abstract}
 We present a new fundamental intuition for why the Kemeny feature of a Markov chain is a constant.  This new perspective has interesting further implications. 
 \end{abstract}
 
\noindent {\bf keywords}: Markov chains, Mixing, Kemeny constant 

\noindent {\bf classification[PACS]}:47A63,35Q70,91G10,11E25

\noindent {\bf classification[MSC]}:60J10, 15B21

\section{Introduction} 
The second-named author has long been interested in the properties of the Kemeny constant in Markov chains, see Hunter \cite{1} and citations therein.  At the 22nd IWMS Conference in Toronto in 2013 he introduced the Kemeny constant to the first-named author and emphasized especially the lack of reasoned, plausible, intuitive argument, apart from purely mathematical justifications, for why this feature of a Markov chain should be a constant.  Subsequently, in Gialampoukidis, Gustafson, Antoniou \cite{2} we accepted its constancy and established the relationship of Kemeny Time to a maximum mixing time for a two-state Markov chain to achieve a total variation distance no greater than any chosen tolerance $ \epsilon $ from the final stationary vector $ \pi $.  Then at the 24th IWMS Conference in Haikou in 2015 the two authors of this paper had further discussions of various issues surrounding the Kemeny constant. As a result of those discussions we found a new intuition from which to view the issue.  The purpose of this short paper is to present that new perspective and some reasoned and plausible supporting arguments. 

The new intuition is to see the well-known basic mean first passage time matrix equation $ M\pi=Ke $ as a change-of-basis procedure.  Once that is carefully written out, but as $ \overline{M}\pi=k $ where we call $ k $ the Kemeny vector, and where $ \overline{M} $ is $ M $ with its diagonal deleted, an insistence on viewing $ \overline{M} $ as the change-of-basis matrix from the $ \overline{M} $ 
column basis to the natural basis, and $ \overline{M}^{{-1}} $ as the change-of-basis matrix from the natural basis to the $ \overline{M} $ column basis, intuits that one must "end up with equally probable pure states". 

For brevity, we will not survey the literature, that having been provided in \cite{1}.  Again for brevity and convenience we will rely upon that paper for notation and basic facts and previously known interpretations of the Kemeny constant in Markov chains.  However, here is some quick background.  
The pioneering book Kemeny and Snell \cite{3} is the origin of the Kemeny feature: the average mean first passage time from any state $ i $ with respect to the equilibrium probability $ \pi $ does not depend on the state $ i $. Here $ P $ is the row-stochastic $ n \times n $ transition matrix for a regular Markov chain with equilibrium (and stationary) probability $ \pi $.  The most relevant pages in \cite{3} are pp. 75-82 and we will refer to those. In particular, the Kemeny feature is embodied in \cite[Theorem~4.4.10]{3}: $ M\alpha^{T}= c\zeta $. This we have written above in more modern notation and as in \cite{1} as $ M\pi = Ke $, where $ e $ is the column vector $ e=(1,\cdots;1)^{T} $ and $ M $ is the matrix $ [m_{ij}] $ of first passage times. $ K $ is commonly called the Kemeny constant and was shown in \cite{3} to be $ K = trace(Z) $ where $ Z = [I-(P-A)]^{{-1}} $ is a resolvent operator and $ A = lim P^{n} $ as $ n\rightarrow \infty $.  

In the ensuing years there arose some disquiet about the meanings of this result and those are detailed in \cite{1}.  A small prize was offered and eventually given to Peter Doyle who showed that the vector components $ k_{i} $ of $ M\pi=k $ satisfy the maximum principle $ k_{i} = \Sigma_{j}p_{ij}k_{j} $ and thus must be constant.  However, this is more the way of proof rather than some deeper intuition so the issue remained still somewhat open. An interesting interpretation of $ K $ as the mean number of links a random surfer will encounter when navigating a random walk on a Markov web until reaching an unknown destination state.  See \cite{1} and \cite{3} for further background information.  

We will prefer to present our new intuition with the always-invertible matrix $ \overline{M} $ which is $ M $ with it's diagonal elements set to zero. This matrix enters also into the proof in \cite{3} and just reduces the Kemeny constant to $ K-1 $.  To conclude this introduction, let us note that it is quite elementary to see from the original treatment in \cite{3} that $ K $ is a constant. From \cite[p. 79]{3} we have $ P(M-D)=M-E $ where $ D $ is the diagonal matrix with elements $ d_{ii}=m_{ii}=\pi_{i}^{{-1}} $ and $ E=ee^{T} $ is the matrix with all the 1's.  Thus $ P\overline{M}=\overline{M}+D-E $ and when applied to $ \pi $ one has
\beq 
P\overline{M}\pi = \overline{M}\pi + e-e(e^{T}\pi ) = \overline{M}\pi \mathrm{.} \label{1.1} 
\eeq
In other words, $ \overline{M}\pi $ is in the principal eigenspace $ sp[e] $ of $ P $ and is therefore a constant times e. 

\section{Why the Kemeny Vector has Equal Coordinates} 

Our approach starts with no Kemeny constant $ K $ at all.  As if we were teaching the introductory linear algebra course, we write the invertible equation $ \overline{M}\pi=k $ as the change of basis:
\begin{equation}
\pi_{1} \begin{bmatrix}
0\\
m_{21}\\
m_{31}
\end{bmatrix} + \pi_{2}
\begin{bmatrix}
m_{12}\\
0\\
m_{32}
\end{bmatrix} + \pi_{3}
\begin{bmatrix}
m_{13}\\
m_{23}\\
0
\end{bmatrix}
= \overline{M}\pi=k=
\begin{bmatrix}
k_{1}\\
k_{2}\\
k_{3}
\end{bmatrix}= k_{1} 
\begin{bmatrix}
1\\
0\\
0
\end{bmatrix} + k_{2}
\begin{bmatrix}
0\\
1\\
0
\end{bmatrix} + k_{3}
\begin{bmatrix}
0\\
0\\
1
\end{bmatrix}
\label{2.1} 
\end{equation}
We have written in three dimensions for clarity but the argument is the same in all dimensions. We call the columns on the left the $ \overline{M} $ column basis and the three columns on he right the natural basis or the pure states or $ e_{1} $, $ e_{2} $, $ e_{3} $ or $ s_{1} $, $ s_{2} $, $ s_{3} $, whatever be your predilection.  


This is why our intuition said: there is an equiprobable pure state assumption somewhere underlying the fact that $ k $ has equal coordinates.  Stated another way, in the way physicists like to claim that one should always work in a "coordinate-free" way: $ \pi $ is "just" $ k $ but now expressed in the $ \overline{M} $ column basis rather than in the pure state "natural" basis.  Stated a third way: the stationary probability $ \pi $, which is the fundamental measure for the process at equilibrium, is really the equiprobability measure in disguise.  

This is a strong claim and a new outcome that we will support in the rest of this paper.  

To begin, our new intuition originated from thinking of (\ref{2.1}) from the change-of-basis procedure as implemented by Gauss row reduction, e.g. see Lay\cite[Section\ 4.7 ]{4}. To invert a matrix equation $ Ax=b $ one forms the tableau $ [A|I] $ and row reduces that to $ [I|A^{{-1}}] $. This is a special case of a general change of basis procedure $ [C|B] \rightarrow [I|\underset{C\leftarrow B}{\operatorname{\textit{P}}}] $ where $ \underset{C\leftarrow B}{\operatorname{\textit{P}}} $ transforms any vector from representation in the $ B $ column basis to its representation in the $ C $ column basis.  In the special case one can say that $ x $ is merely $ b $ changed from its representation in the natural basis to it's representation in the $ A $ column basis.  

We will illustrate this in the next section by explicitly carrying it out for the Land of Oz example of \cite{3}.  

Of course the change-of-basis matrix inversion perspective applied to $ \overline{M}\pi=k $ and the $ \pi=\overline{M}^{{-1}}k $ is just a special case of representing any vector $ b $ written as usual in the natural basis to changing it's representation to $ x=A^{{-1}}b $ where $ x $ is now its coordinates in the $ A $ column basis. The key here is that $ \pi $ is a very special equilibrium probability measure. 

\section{The Change-of-Basis Picture} 
Because our new intuition arose out of insisting that we view the remarkable Kemeny-Snell equation $ M\pi=Ke $ as a change-of-basis statement, we elaborate by specific example here.  A good elementary reference is the book \cite[Section\ 4.7\ pp\ 239-242]{4}.  We may immediately get into the spirit by doing the key example used throughout \cite{3}: the Land of Oz example
\begin{equation}P=
\begin{array}{r}
\begin{matrix}
R & N & S
\end{matrix} \\
\begin{matrix}
R\\
N\\
S
\end{matrix} 
 \begin{bmatrix}
\frac{1}{2} & \frac{1}{4} & \frac{1}{4}\\
\frac{1}{2} & 0 & \frac{1}{2}\\
\frac{1}{4} & \frac{1}{4} & \frac{1}{2}\\
\end{bmatrix}\ 
\end{array}
\label{3.1} 
\end{equation}
We know that $ Pe=e $, $ P^{T}\pi=\pi=(\pi_{1},\pi_{2}, \pi_{3})^{T}=(\frac{2}{5},\frac{1}{5},\frac{2}{5})^{T} $, and as computed in \cite{3} via the resolvent operator $ Z $, the mean first passage time matrix $ M $ is

\begin{equation}
M= \begin{bmatrix}
\frac{5}{2} & 4 &\frac{10}{3}\\
\frac{8}{3} & 5 & \frac{8}{3}\\
\frac{10}{3} & 4 &\frac{5}{2}
\end{bmatrix} 
\label{3.2}  
\end{equation}
To calculate $ \overline{M}^{{{-1}}} $ by the Gauss procedure, one row reduces the tableau as follows:
\begin{equation}
\begin{array}{c}
\left[ \begin{array} {@{}ccc|ccc@{}}
0 & 4 & \frac{10}{3} & 1 & 0 & 0 \\
\frac{8}{3} & 0 & \frac{8}{3} & 0 & 1 & 0 \\
\frac{10}{3} & 4 & 0 & 0 & 0 &1 \\
\end{array}
\right] \rightarrow
\\
\begin{matrix}
M_{1}\, M_{2}\,M_{3}
&
s_{1}\,s_{2}\,s_{3}
\end{matrix}
\end{array}
\begin{array}{c}
\left[ \begin{array} {@{}ccc|ccc@{}}
1 & 0 & 0 & \frac{{-3}}{20} & \frac{3}{16}& \frac{3}{20} \\
0 & 1 & 0 & \frac{1}{8} & \frac{{-5}}{32} & \frac{1}{8} \\
0 & 0 & 1 & \frac{3}{20} & \frac{3}{16} & \frac{{-3}}{20} \\
\end{array}
\right]\\
\begin{matrix}
\, \, 
\end{matrix}
\end{array} 
\label{3.3}
\end{equation}
This is a special case of the more general change-of-basis in which one drives the tableau:\\
$ [C_{1}\  C_{2}\   \cdots\   C_{n}\  | B_{1}\ B_{2}\ \cdots\  B_{n}] \rightarrow [I|\underset{C\leftarrow B}{\operatorname{\textit{P}}}]$.  

Thus the Land of Oz Markov chain Mean first passage time matrix $ \overline{M}^{{{-1}}} $on the right side of (\ref{3.3}) exactly changes the representation of vectors in the natural basis $ \{s_{1},s_{2},s_{3}\} $ of pure states into representations in terms of the mean first passage time column basis $ \{\overline{M}_{1}, \overline{M}_{2}, \overline{M}_{3}\} $. In particular, $ \overline{M}^{{-1}} $ transforms the equally probable measure $ \frac{e}{3}=(\frac{1}{3},\frac{1}{3},\frac{1}{3})^{T} $ to a multiple of the stationary measure $ (\pi_{1},\pi_{2},\pi_{3}) $.  Generally for the $ n \times n $ case where $ \overline{M}\pi= (K-1)e $, we make the right side of measure one by dividing both sides by $ n(K-1) $ a factor which can be absorbed by $ \overline{M} $ and its inverse.  One easily calculates that $ K-1 =\frac{32}{15} $ for the Land of Oz chain so the normalizing factor is $ \frac{32}{5} $.  

While this change-of-basis picture brings to the fore that the right side of (\ref{2.1}) is actually a representation of the Kemeny-Snell vector $ k $ in terms of the pure states $ s_{1}, s_{2}, s_{3} $, it does not prove that $ k_{1}=k_{2}=k_{3} $.  That fact was already established in \cite{3} and has been shown other ways, see \cite{1}.  We gave a very simple proof at the end of Section 1.  Here is another one, which we wish to mention in order to bring us to the point we emphasized at the end of Section 2: $ \pi $ is a very special vector measure-theoretically.

Just apply $ P^{n} $ to both sides of the change-of-basis equation (\ref{2.1}) and go the limit as $ n \rightarrow \infty $.  The left side is invariant since  $ P^{n}(\overline{M}\pi)=\overline{M}\pi $ as we showed in Section 1.  The right converges to
\begin{equation}
k_{1}e\pi^{T} \begin{bmatrix}
1\\
0\\
0
\end{bmatrix} + k_{2}e\pi^{T} \begin{bmatrix}
0\\
1\\
0
\end{bmatrix} + k_{3}e\pi^{T} \begin{bmatrix}
0\\
0\\
1
\end{bmatrix}
= (k_{1}\pi_{1}+k_{2}\pi_{2}+k_{3}\pi_{3})e 
\label{3.4}
\end{equation} so that left side $ \overline{M}\pi=k $ is a constant multiple $ (K-1) $ of $ e $.  Here we have used the fact that $ \lim P^{n}=A $ in the Kemeny-Snell notation \cite{3} is the rank-one oblique projection given by $ A=e\pi^{T} $.

For the Perron convergence theory see Hunter \cite[Chapter\ 7]{5} and Horn and Johnson \cite[Chapter\ 8]{6} and especially their wonderful Lemma 8.2.7 on pages 497-498.  In their notation $ \lim P^{n} $ is $ L=xy^{T}= e\pi^{T} $ here and we sometimes like to go further, see Gustafson \cite[p.\ 206]{7} to regard the normalized version $ \frac{xy^{T}}{y^{T}x} $ as the oblique projection onto the span of $ x $ from the direction perpendicular to $ y $.  That $ L^{2}=L $ projection view provides a strictly geometrical new view of $ K $: the amplitude of the oblique rank-one projection $ L(M\pi) $ onto $ sp[e] $.  

Thus the change-of-basis equation (\ref{2.1}) by the invariance of its left side $ \overline{M}\pi $ under the Markov chain's transition matrix iterates $ P^{m} $ as shown in equation (\ref{1.1}) has led us to the fact (\ref{3.3}) that the Markov process must "end up with equally probable pure states".  The later are the essence of the at-first seemingly harmless eigenvector $ e $.  The fact this occurs rests principally upon the stationary probability $ \pi $.  
\section{Discussion} 
Our new perspective raises a number of interesting implications.  Some of these may be worthy of further study but we can only mention a couple of them here in this brief paper.

Why equi-probability? The reply: Kemeny-Snell's \cite{3} remarkable equation $ M\pi=Ke $ is only a statement at equilibrium.  Everyone knows that one can start a regular Markov chain with any initial probability and iterate until you get to the limit distribution $ \pi $.  This is generalized in the famous Perron Theorem, e.g, \cite[p.\ 499]{6}, and the point is that the $ L^{\infty} $ limit of $ P^{n}$ is $ L=xy^{T}=e\pi^{T} $ in our case. $ L $ is a rank-one oblique projector and in fact it itself represents an independent trials process with transition matrix
\begin{equation}
L=e\pi^{T}= \begin{bmatrix}
\pi_{1} & \pi_{2} & \pi_{3}\\
\pi_{1} & \pi_{2} & \pi_{3}\\
\pi_{1} & \pi_{2} & \pi_{3}\\
\end{bmatrix}
\end{equation}\label{4.1}
with Perron eigenvector $ Le=e $  and stationary equilibrium probability $ L^{T}\pi=\pi e^{T}\pi = \pi $.  

An MCMC implication? The widely acclaimed Markov Chain Monte Carlo, see e.g Antoniou, Christidis, Gustafson \cite{8}, assumes you can find an initial distribution $ \pi_{0} $ which after a sufficient number of interations is close to the invariant distribution $ \pi $ which is believed to represent the physical process being modeled.  One then performs Monte Carlo simulations on the latter.  Our interpretation in \cite{8} is that the iterations generate sufficient mixing so that the subsequent sampling stage represents adequately the regular probability distribution of the application.  We go further \cite{8} and hope that there exists a deeper underlying physical dynamics. Here we say: do your Monte Carlo equiprobably.   

Next, we mention that we became curious about how Kemeny-Snell \cite{3} somehow were able to move effortlessly between $ P $ and $ P^{T} $, or if you wish between $ M $ and $ M^{T} $, $ viz $, between \cite[Theorems\ 4.4.9\ and\ 4.4.10]{3}.  The technical secret seems to lie in the second term in equation (\ref{1.1}) in our introduction.  Namely, the symetric operator $ D-E $ has null space $ sp\{\pi\} $.  One could go a bit further intuitively and assert that $ D $ represents the probability of the self loops of the pure state $ s_{1},s_{2},s_{3} $ and $ E $ represents random equiprobable noise and the two are canceled on the stationary distribution $ \pi $.  

We may ask how our column bases (the columns of $ M $) behave as the Markov process progresses.  That is, we expect Kemeny time $ K $ to 'decrease' as we step forward in the chain $ P,P^{2},\cdots $.  To make this precise, recall $ K=1+\sum_{2}^{n}(1-\lambda_{i})^{{-1}} $, and let us make the additional assumption that $ P $ is primitive so that all the $ | \lambda_{i} | < 1 $ for all of $ i>1 $. The Kemeny-Snell equation $ M_{m}\pi=k_{m}=K_{m}e $ at the $ m $th step in the Chain has Kemeny time $ K_{m}=1+\sum^{n}_{i=2}(1-\lambda^{m}_{i})^{{-1}} $ which converges down to $ K_{L}=n $ as the $ |\lambda_{i}|^{m} $ all go to zero.  The column bases of $ M_{m} $ converge to those of $ M_{L} $ which for $ n=3 $ are
\begin{equation}
M_{L}=
\begin{bmatrix}
\pi_{1}^{{-1}} \begin{bmatrix}
1\\
1\\
1
\end{bmatrix} &
\pi_{2}^{{-1}} \begin{bmatrix}
1\\
1\\
1
\end{bmatrix} &
\pi_{3}^{{-1}} \begin{bmatrix}
1\\
1\\
1
\end{bmatrix}
\end{bmatrix} = e [ \pi_{1}^{{-1}}, \pi_{2}^{{-1}}, \pi_{3}^{{-1}} ]
\label{4.2}
\end{equation} Notice that for the Land of Oz examples (see Section 3) this means that some of the mean first passage times $ m_{ij} $ increase while others decrease as the $ M_{m} $ converge toward $ M_{L}=e[\frac{5}{2},5,\frac{5}{2}] $.  The latter is a rank-one matrix, so its columns are no longer a basis even if those of the $ M_{m} $ were, but there is no problem with $ \overline{M}_{L} $ which conserves our change of basis picture $ \overline{M}_{L}\pi=2e $ in this and all examples.  
\section{Conclusions} 
In the recent paper \cite{1} and before that it has been emphasized that there was still needed a better reasoned, plausible intuitive argument, apart from purely mathematical justifications, for why the Kemeny feature of a Markov chain should be constant.  Here we have shared with you a new intuition, reasoned arguments supporting that intuition, and a perhaps unexpected plausible fundamental outcome.  The intuition was to insist on viewing the remarkable Kemeny-Snell first passage time equation $ M\pi=k $ as an $ M $-column basis representation of $ k $, then wonder why the new coordinates $ k_{1},k_{2},k_{3} $ of the natural basis representation of $ \pi $ need to be equal.  Of course that perspective holds for arbitrary dimension $ n $.  The resulting reasoned arguments followed closely the original treatment in \cite{3} and, by the way, completely avoided the machineries of operator resolvents or generalized group inverses.  The other perspective in our reasoned arguments was the Perron Theorem and especially its limit oblique projection $ e\pi^{T} $.  The plausible outcome was that the Markov chain in the limit must converge to equally probable pure states.  This equiprobability measure is hidden within the equilibrium measure $ \pi $.  In important applications it is postulated to represent a deeper underlying chaos \cite{8}.

\section*{Acknowledgement}
Karl Gustafson would like to express his thanks to Jeff Hunter and Simo Puntanen as chairs of the IWMS-2015 for inviting him and to ILAS for designating him as their lecturer for the conference.

\end{document}